\documentclass[11pt]{article}
\usepackage{amssymb,amsmath,mathrsfs,latexsym,amscd}
\usepackage{exscale,latexsym,amsthm,graphics}

\usepackage{makeidx}
\makeindex
\usepackage{epsfig}

\newtheorem{defn}{Definition}[section]
\newtheorem{thm}[defn]{THEOREM}
\newtheorem{lemma}[defn]{LEMMA}
\newtheorem{prop}[defn]{PROPOSITION}
\newtheorem{proposition}[defn]{PROPOSITION}

\newtheorem{remark}[defn]{{\it Remark}}

\usepackage{txfonts}

\def\pf{{\medskip\noindent { Proof. }}}
\def\qed{{\hfill $\Box$ \bigskip}}

\def\bM {{\mathbb M}}

\def\N {{\mathbb N}}
\def\R {{\mathbb R}}

\def\EE{{\mathbb E}}
\def\P{{\mathbb P}}

\def\E{{\mathcal E}}
\def\P{{\mathbb P}}

\def\grad{{\nabla}}

\begin{document}

\noindent
\centerline{\Large\bf POINTWISE WEAK EXISTENCE OF DISTORTED }
\centerline{\Large\bf SKEW BROWNIAN MOTION WITH RESPECT TO}
\centerline{\Large\bf DISCONTINUOUS MUCKENHOUPT WEIGHTS}
\bigskip
\centerline{}
\centerline{ JIYONG SHIN and GERALD TRUTNAU} 

\bigskip
\noindent

\bigskip
\noindent
{For any starting point in $\mathbb{R}^d$, $d \ge 3$, we identify the stochastic differential equation of distorted Brownian motion with respect to a certain discontinuous Muckenhoupt $A_2$-weight $\psi$. 
The discontinuities of $\psi$ typically take place on a sequence of level sets of the Euclidean norm 
$D_k:=\{x\in \mathbb{R}^d\, | \;\|x\|=d_k\}$, $k\in\mathbb{Z}$, where $(d_k)_{k\in\mathbb{Z}}\subset (0,\infty)$ 
may have accumulation points and each level set $D_k$ plays the role of a permeable membrane.}



\section{INTRODUCTION}
In this note we are concerned with the construction of a weak solution to the stochastic differential equation 
\begin{eqnarray}\label{eq1}
X_t=x+W_t+\int_0^t \frac{\nabla \rho}{2\rho}(X_s) \ ds+\int_0^{\infty}\int_0^t \nu_{a}(X_s) \, d\ell_s^a(\|X\|)\,\eta(da),\ \ x\in \mathbb{R}^d,
\end{eqnarray}
where $W$ is a $d$-dimensional standard Brownian motion, $\rho$ is typically a Muckenhoupt $A_2$-weight, $\nu_{a}$ 
is the unit outward normal on the boundary of the Euclidean ball of radius $a$ about zero, $\ell^{a}(\|X\|)$ is the 
symmetric semimartingale local time at $a\in (0,\infty)$ of $\|X\|$, $\eta=\sum_{k\in \mathbb{Z}}(2\alpha_k-1)\delta_{d_k}$ with $(\alpha_k)_{k\in \mathbb{Z}}\subset (0,1)$ is a sum of 
Dirac measures at a sequence  $(d_k)_{k\in \mathbb{Z}}\subset (0,\infty)$  
with exactly two accumulation points in $[0,\infty)$, one is zero and the other is $m_0>0$. More accumulation points 
could be allowed. For a discussion on this point we refer to \cite[Remark 2.6(ii)]{OuRuTr13}. The absolutely continuous 
component of drift $\frac{\nabla \rho}{2\rho}$ is typically unbounded and discontinuous. For an interpretation of the equation 
we refer to Theorem \ref{T;MAIN} and Remark \ref{interpretation}. Variants of (\ref{eq1}) with reflection on hyperplanes, 
instead of balls, but without accumulation points and Lipschitz drift appear in \cite{Tom, O81, O82, taka86a}.
In particular, (\ref{eq1}) is a multidimensional analogue of an equation that was thoroughly studied in \cite{OuRuTr13} 
and both equations share a lot of similarities. For instance the 
way to determine $(\gamma_{k})_{k\in \mathbb{Z}}$ and $(\overline{\gamma}_{k})_{k\in \mathbb{Z}}$ in (\ref{eqphi}) below, when 
$(\alpha_k)_{k\in \mathbb{Z}}$ in (\ref{eq1}) is given, is the same here as in \cite[Proposition 2.11]{OuRuTr13}. 
The way to obtain  $(\alpha_k)_{k\in \mathbb{Z}}$ from $(\gamma_{k})_{k\in \mathbb{Z}}$ and $(\overline{\gamma}_{k})_{k\in \mathbb{Z}}$ 
is described in Remark \ref{interpretation} (cf. also \cite[Remark 2.4(ii) and Theorem 2.5]{OuRuTr13}). See further \cite[Remark 2.7]{OuRuTr13} 
for another similarity and \cite[Section 3.3]{OuRuTr13} as well as references therein for a possible application to models with countably many permeable membranes that accumulate. 
For the construction of a solution to (\ref{eq1}) for any starting 
point $x\in \mathbb{R}^d$ the key point is to identify (\ref{eq1}) as distorted Brownian motion (see \cite{ahkrstr}, \cite{fuku81} 
for an introduction to distorted Brownian motion). Then, one needs to show that the absolute continuity condition \cite[p.165]{FOT} is 
satisfied for the underlying Dirichlet form and that the strict Fukushima decomposition \cite[Theorem 5.5.5]{FOT} is applicable. 
In order to identify (\ref{eq1}) as distorted Brownian motion we proceed informally as follows. 
We consider the density $\rho\phi$ for the underlying Dirichlet form in (\ref{DF}), where $\phi$ is a step function on annuli 
in $\mathbb{R}^d$, see (\ref{eqphi}), and $\rho$ is typically a Muckenhoupt $A_2$-weight (see Remark \ref{R;RE1}(iii) below). 
For the precise conditions we refer to (H1)-(H3) in section \ref{2}. 
The logarithmic derivative, which is the drift of the distorted Brownian motion, is then informally given by
$$
\frac{d(\rho\phi)}{\rho\phi}=\frac{d\rho}{\rho}+\frac{d\phi}{\phi}
$$
(cf. \cite[Remark 2.6]{OuRuTr13}). This is rigorously performed through an integration by parts formula in Proposition \ref{ibp} below. 
By results on Muckenhoupt weights in \cite{St2,St3}, the existence of a jointly 
continuous transition kernel density for the semigroup associated to the Dirichlet form given in 
(\ref{DF}) is obtained. A Hunt process with the given transition kernel density is implicitly assumed to exist in condition (H4) of section \ref{2} (for ways to construct such a process, we refer to \cite{ShTr}, see also Remark \ref{remarkAKR}). 
Under the conditions (H1)-(H4), we then show in a series of statements in sections \ref{4} and \ref{5}, that the strict Fukushima decomposition can be applied to obtain a solution to (\ref{eq1}) (see main Theorem \ref{T;MAIN}). \\
Finally one could think of generalizing (\ref{eq1}), or more precisely the process of Theorem \ref{T;MAIN}(i) with reflections on boundaries of Lipschitz domains (instead of smooth Euclidean balls). The main ingredient to obtain this generalization would be \cite[Theorem 5.1]{Tr1} (see \cite[section 5]{Tr3}). 
In case of skew reflection at the boundary of a single $C^{1,\lambda}$-domain, $\lambda\in(0,1]$, $\rho\equiv 1$, and smooth diffusion coefficient, a weak solution has been constructed in \cite[III. \S 3 and \S 4]{port}, see also references therein. However, the reflection term is defined as generalized drift and not explicitly as in Theorem \ref{T;MAIN}.
\section{FORMULATION OF THE MAIN THEOREM}\label{2}
Let $m_0 \in (0, \infty)$ and $(l_{k})_{k \in \mathbb{Z}} \subset (0, m_0)$, $0 < l_{k} < l_{k+1} <m_0$, be a sequence converging to $0$ as $k \rightarrow -\infty$ and converging to $m_0$ as $k \rightarrow \infty$. Let $(r_{k})_{k \in \mathbb{Z}} \subset (m_0, \infty)$, $m_0 < r_{k} < r_{k+1} < \infty$, be a sequence converging to $m_0$ as $k \rightarrow -\infty$ and tending to infinity  as $k \rightarrow \infty$. Let
\begin{eqnarray}\label{eqphi}
\phi : = \sum_{k \in \mathbb{Z} } \left (\gamma_{k}  \cdot 1_{A_{k}}
+ \overline{\gamma}_{k} \cdot 1_{\hat{A}_{k}}\right ),
\end{eqnarray}
where $\gamma_{k}$ , $ \overline{\gamma}_{k} \in (0, \infty)$, $A_{k} : = B_{l_{k}} \setminus \overline{B}_{l_{k-1}}$ , $\hat{A}_{k} : = B_{r_{k}} \setminus \overline{B}_{r_{k-1}}$ , $B_{r}: = \{ x \in\R^d \,|\; \|x\| < r  \}$ for any $r>0$, $\overline{B}_r$ denotes the closure of $B_{r}$, $1_{A}$ the indicator function of a set $A$, and $\|x\|$ the Euclidean norm of $x\in \mathbb{R}^d$. We denote by $d\sigma_r$ the surface measure on the boundary $\partial B_r$ of $B_r$, $r>0$. \\
Let $d\ge 3$, and let $C_0^{\infty}(\R^d)$ denote the set of all infinitely differentiable functions with compact support in $\R^d$. Let $\nabla f : = ( \partial_{1} f, \dots , \partial_{d} f )$ where $\partial_j f$ is the $j$-th partial derivative of $f$, and $\Delta f : = \sum_{j=1}^{d} \partial_{jj} f$ and $\partial_{jj} f := \partial_{j}(\partial_{j} f) $, $j=1, \dots, d$. 
As usual $dx$ is the Lebesgue measure on $\R^d$ and  $L^q(\R^d, \mu)$, $q \ge 1$, are the usual $L^q$-spaces with respect to the  measure $\mu$ on $\R^d$, and $L^{q}_{loc}(\R^d,\mu) := \{ f \,|\; f \cdot 1_{U} \in 
L^q(\R^d, \mu),\,\forall U \subset \R^d, U \text{ relatively compact open}   \}$. 
For any open set $\Omega \subset \R^d$, $H^{1,q}(\Omega, dx), \; q \ge 1$ is defined to be the set of all functions $f \in L^{q}(\Omega, dx)$ such that $\partial_{j} f \in L^{q}(\Omega, dx)$, $j=1, \dots, d$, and 
$$
H^{1,q}_{loc}(\R^d, dx) : =  \{ f  \,|\;  f \cdot 1_{U} \in H^{1,q}(U, dx),\,\forall U \subset \R^d,\, U \text{ relatively compact open}  \}.
$$
For a topological space X $\subset \R^d$ with Borel  $\sigma$-algebra $\mathcal{B}(X)$ we denote the set of all $\mathcal{B}(X)$-measurable $f : X \rightarrow \R$ which are bounded, or nonnegative by $\mathcal{B}_b(X)$, $\mathcal{B}^{+}(X)$ respectively. We further denote the set of continuous functions on $X$, the set of continuous bounded functions on $X$ by $C(X)$, $C_b(X)$ respectively.\\ 
A function $\psi\in \mathcal{B}(\R^d)$ with $\psi>0$ $dx$-a.e. is said to be a Muckenhoupt $A_2$-weight (in notation $\psi\in A_2$), if there exists a positive constant $A$ such that, for every ball $B\subset\R^d$,
\[
\left(\int_B \psi dx\right) \left(\int_B \psi^{-1} dx\right) \le A \, \left(  \int_{B} 1\, dx \right)^2.
\]
For more on Muckenhoupt weights, we refer to \cite{Tu}. \\
Throughout we shall assume
\begin{itemize}
\item[(H1)] $\sum_{k \in \mathbb{Z}} |\,  \gamma_{k+1} - \gamma_{k} \,| + \sum_{k \le 0} |\,  \overline{\gamma}_{k+1} - \overline{\gamma}_{k} \,| < \infty$
and for all $r>0$ there exists $\delta_{r}>0$  such that $\phi \ge \delta_{r}$ $dx$-a.e. on $B_{r}$.
\item[(H2)] $\rho \, \phi \in A_2$.
\end{itemize}

\begin{remark}\label{R;RE1}
(i) (H1) implies that $\phi \in L^{1}_{loc}(\R^d, dx)$ and that $\gamma : = \lim_{k \rightarrow \infty} \gamma_{k}$, $\overline{\gamma} : = \lim_{k \rightarrow -\infty} \overline{\gamma}_{k}$  exist and  $\gamma>0$, $\overline{\gamma}>0$. In particular,  $\phi$ is locally bounded above and locally bounded away from zero.\\
(ii) (H1) and (H2) imply $\rho>0$ $dx$-a.e.\\
(iii) Let $c>0$. If $c^{-1} \le \phi \le c$ and $\rho \in A_{2}$, or if $\rho = 1$ and  $c^{-1} \, \psi \le \phi \le c \, \psi$ for some $\psi \in A_{2}$, then $\rho \, \phi \in A_{2}$.
\end{remark}
Furthermore, we shall throughout assume the following condition.
\begin{itemize}
\item[(H3)] $\rho = \xi^2$ for some $ \xi \in H^{1,2}_{loc}(\R^d, dx)$.
\end{itemize}
\begin{remark}\label{rem2}
(H3) implies that $\rho \in H^{1,1}_{loc}(\R^d, dx)$ and by (H1) $\frac{\|\nabla \rho \|}{\rho} \in L^{2}_{loc}( \R^{d} , \; \rho \phi \, dx)$ since $\phi$ is locally bounded above $dx$-a.e.
\end{remark}
We consider the symmetric positive definite bilinear form
\begin{equation}\label{DF}
\E (f,g) : = \frac{1}{2}\int_{\R^d}\nabla f\cdot\nabla g\,  \left( \rho \, \phi \right) \,dx , \ \ f,g \in C_0^{\infty}(\R^d).
\end{equation}
Since $\rho \, \phi \in A_2$, we have $\frac{1}{\rho \, \phi}  \in L^1_{loc}(\R^d, dx)$, and the latter implies that  $\eqref{DF}$  is closable in $L^2(\R^d, \rho  \phi \,dx)$ (see \cite[II.2 a)]{MR}). The closure $(\E,D(\E))$ of $\eqref{DF}$ is a strongly local, regular, symmetric Dirichlet form (cf. e.g. \cite[p. 274]{St3}). As usual we define $\E_1(f,g) := \E(f,g) + (f,g)_{L^{2}(\R^d , \, \rho \phi \, dx)}$ for $f,g \in D(\E)$ and  $\| \,f\, \|_{D(\mathcal{E})} : = \E_1(f,f)^{1/2},  \;   f \in D(\E)$.
Let  $(T_t)_{t \ge 0}$ and $(G_{\alpha})_{\alpha > 0}$ be the $L^2(\R^d, \rho \phi \, dx)$-semigroup and resolvent associated to $(\E,D(\E))$ and $(L,D(L))$ be the corresponding generator (see \cite[Diagram 3, p. 39]{MR}).
From \cite[p. 303 Proposition 2.3]{St2} and \cite[p. 286 A)]{St3} we know that there exists a jointly continuous transition kernel density  $p_{t}(x,y)$  such that
\[
P_t f(x) := \int_{\R^d} p_t(x,y)\, f(y)\;\rho(y) \, \phi(y) \,dy \,, \;\; t>0, \;\;x \in \R^d, 
\]
$f\in \mathcal{B}_b(\R^d)$, is a $\rho \phi dy$-version of $T_t f$ if $f  \in  L^2(\R^d , \rho \, \phi \, dx) \cap \mathcal{B}_b(\R^d)$. Furthermore, taking the Laplace transform of $p_t(x,y)$ there exists a resolvent kernel density $r_{\alpha}(x,y)$ such that
\[
R_{\alpha} f(x) := \int_{\R^d} r_{\alpha}(x,y)\, f(y) \;\rho(y) \, \phi(y) \,dy \,, \;\; \alpha>0, \;\;x \in \R^d, 
\]
$f\in \mathcal{B}_b(\R^d)$, is a $\rho \phi dy$-version of $G_{\alpha} f$ if $f  \in  L^2(\R^d , \rho \, \phi \, dx) \cap \mathcal{B}_b(\R^d)$. Accordingly, for a signed Radon measure $\mu$,
let us define
\[
R_{\alpha} \mu (x) = \int_{\R^d} r_{\alpha}(x,y) \, \mu(dy) \, , \;\; \alpha>0, \;\;x \in \R^d,
\]
whenever this makes sense. Since $p_t(\cdot, \cdot)$ is jointly continuous and $p_t(x,y)$ has exponential decay in $y$ for fixed $t$ and $x$ in a compact set, $(P_t)_{t \ge 0}$ is strong Feller, i.e. $P_t(\mathcal{B}_b(\R^d)) \subset C_b(\R^d)$. For details, we refer to \cite{ShTr}. In particular $R_1(\mathcal{B}_b(\R^d)) \subset C_b(\R^d)$.\\
We consider further the following condition
\begin{itemize}
\item[(H4)] There exists a Hunt process \begin{equation*}
\bM = (\mathbf{\Omega} , \mathcal{F}, (\mathcal{F}_t)_{t\geq0}, \zeta ,(X_t)_{t\geq0} , (\P_x)_{x \in \R^d \cup \{\Delta\} }),
\end{equation*}
with state space $\R^d \cup \{\Delta\}$ and life time $\zeta$, which has $(P_t)_{t \ge 0}$ as transition function, $R_1f$ 
is continuous for any $f \in L^2(\R^d, \rho \phi dx)$ with compact support, and if $\phi \nequiv$ const., 
we additionally assume $R_1(1_G \ \rho d \sigma_r)$ is continuous for any $G \subset \R^d$ relatively compact open, $r > 0$.
\end{itemize}
In (H4), $\Delta$ is the cemetery 
point and as usual any function $f : \R^d \rightarrow \R$ is extended to $\{\Delta\}$ by setting $f(\Delta):=0$. 

\begin{remark}\label{remarkAKR} 
A Hunt process associated with $(P_t)_{t \ge 0}$ as in (H4) can be constructed by the same methods as used in \cite[section 4]{AKR}. 
The method used in \cite[section 4]{AKR} is applicable, if one can find enough nice functions in  $D(L)$ (cf. \cite[proof of Lemma 4.6]{AKR}). A Hunt process with transition function $(P_t)_{t \ge 0}$ as in (H4) can also be constructed by showing that $(P_t)_{t \ge 0}$ defines a classical Feller semigroup. For details and concrete examples we refer to \cite{ShTr}.
\end{remark}

Under (H4), $\bM$ satisfies in particular the {\it absolute continuity condition} as stated in \cite[p. 165]{FOT}. 

\begin{proposition}\label{p;conserv}
Let a Dirichlet form be given as the closure of
\begin{equation*}
\frac{1}{2}\int_{\R^d}\nabla f\cdot\nabla g  \;  \psi  \,dx , \ \ f,g \in C_0^{\infty}(\R^d)
\end{equation*}
in $L^2(\R^d, \psi \, dx)$ where $\psi \in A_{2}$. Then it is conservative.
\end{proposition}
\pf
By \cite[Proposition 1.2.7]{Tu} $\psi \, dx$ is volume doubling. Hence by \cite[Proposition 5.1, Proposition 5.2]{GrHu}
$$
c_{1} \, r^{\alpha'} \le \psi \, dx(B_{r}) \le c_{2} \, r^{\alpha}, \;\; \forall r \ge 1,
$$
where $c_{1}, c_{2}, \alpha, \alpha' > 0$ are constants. In particular,
$$
\int_{1}^{\infty} \frac{r}{\log \big ( \psi \, dx(B_{r})  \big )} \, dr =\infty.
$$
Hence conservativeness follows by \cite[Theorem 4]{St1}.
\qed

\begin{remark}
Proposition $\ref{p;conserv}$ holds more generally for $A_{p}$-weights, $p \in [1,\infty)$ (see \cite[Definitoin 1.2.2]{Tu} for the definition of $A_{p}$) by the same arguments as in the proof of Proposition $\ref{p;conserv}$.
\end{remark}
It follows from  Proposition \ref{p;conserv} and the strong Feller property of $(P_t)_{t \ge 0}$ that under (H4)
\begin{equation}\label{HPcons}
\P_{x}(\zeta = \infty)=1, \ \ \forall x \in \R^d.
\end{equation} 
We will refer to \cite{FOT} from now on till the end, hence some of its standard notations may be adopted below without definition. The following theorem is the main result of our paper.
It will be proved in section \ref{5}.
\begin{thm}\label{T;MAIN}
Suppose (H1)-(H4) hold. Then: \\
(i) The process $\bM$ satisfies
\begin{equation}\label{SDE}
X_t = x + W_t + \int^{t}_{0} \frac{\grad{\rho}}{2\rho} (X_s) \, ds +  N_{t} \; , \;\; t\ge 0,
\end{equation}
$ \P_x $ -a.s. for any $x \in \R^d$, where $W$ is a standard d-dimensional Brownian motion starting from zero and
\begin{equation*}
N_t : =  \sum_{k \in \mathbb{Z}} \left ( \frac{ \gamma_{k+1} - \gamma_{k}  } {  \gamma_{k+1} + \gamma_{k}}  \int_{0}^{t} \nu_{l_{k}}(X_s) \; d \ell_s^{l_{k}}
+ \frac{  \overline{\gamma}_{k+1} - \overline{\gamma}_{k} }{  \overline{\gamma}_{k+1} + \overline{\gamma}_{k} } \ \int_{0}^{t} \nu_{r_{k}}(X_s) \; d \ell_s^{r_{k}}\right )
\end{equation*}
\begin{equation*}
+ \; \frac{ \overline{\gamma}  - \gamma }{ \overline{\gamma}  + \gamma } \int_{0}^{t} \nu_{m_0} (X_s) \; d \ell_s^{m_0},
\end{equation*}
where $\nu_{r} = (\nu_{r}^1 , \dots , \nu_{r}^d)$, $r>0$ is the unit outward normal vector on $\partial  B_{r}$ and $\ell^{l_{k}}$, $\ell^{r_{k}}$ and $\ell^{m_0}$ are boundary local times of $X$, i.e. they are positive continuous additive functionals of $X$  in the strict sense associated via the Revuz correspondence (cf. \cite[Theorem 5.1.3]{FOT}) with the weighted surface measures $\frac{\gamma_{k+1} + \gamma_{k}}{2} \,  \rho \, d \sigma_{l_{k}}$ on $\partial B_{l_{k}}$,  $\frac{\overline{\gamma}_{k+1} +  \overline{\gamma}_{k} }{2} \,  \rho \, d \sigma_{r_{k}}$ on $\partial B_{r_{k}}$, and  $\frac{\overline{\gamma} + \gamma}{2} \,  \rho \, d \sigma_{m_0}$ on $\partial B_{m_0}$, respectively, and related via the formulas
\begin{eqnarray*}
\EE_x \left [ \int_{0}^{\infty} e^{-t} \, d \ell_{t}^{l_k}   \right] &=& R_1 \left( \frac{\gamma_{k+1} + \gamma_{k}}{2} \, \rho d\sigma_{l_{k}} \right)(x),\\
\EE_x \left[ \int_{0}^{\infty} e^{-t} \, d \ell_{t}^{r_k}   \right] &=& R_1 \left( \frac{\overline{\gamma}_{k+1} +  \overline{\gamma}_{k} }{2}  \, \rho d\sigma_{r_{k}} \right)(x),\\
\EE_x \left[ \int_{0}^{\infty} e^{-t} \, d \ell_{t}^{m_0}   \right] &=& R_1 \left( \frac{\overline{\gamma} + \gamma}{2} \, \rho d\sigma_{m_0} \right)(x),
\end{eqnarray*}
which all hold for any $x \in \R^d$, $k \in \mathbb{Z}$.\\
(ii) $\big ( (\|X_t\|)_{t\ge 0}, \P_x \big )$ is a continuous semimartingale for any $x\in \mathbb{R}^d$ and 
$$
\P_x \big ( \ell_{t}^{a}= \ell_{t}^{a}(\|X\|) \big )=1,\ \ \ \forall x\in \mathbb{R}^d, \ t\ge 0, \ a\in\{m_0, l_k,r_k : k\in \mathbb{Z}\},
$$
where $\ell_{t}^{a}(\|X\|)$ is the symmetric semimartingale local time of $\|X\|$ at $a\in (0,\infty)$ as defined in \cite[VI.(1.25)]{RYor}.
\end{thm}

\begin{remark}\label{interpretation}
In view of Theorem \ref{T;MAIN}, the non absolutely continuous drift component $N$ in (\ref{SDE}) may be interpreted as follows. Define
\[
\alpha_{k}:= \frac{\gamma_{k+1}}{\gamma_{k+1} + \gamma_{k}}, \ \ \ k\in\mathbb{Z}.
\]
Then $\alpha_{k}\in (0,1)$ and 
$$
\frac{ \gamma_{k+1} - \gamma_{k} } { \gamma_{k+1} + \gamma_{k} }=2\alpha_{k}-1=\alpha_{k}-(1-\alpha_{k}).
$$
and so the drift component 
$$
\frac{ \gamma_{k+1} - \gamma_{k} } {  \gamma_{k+1} + \gamma_{k} }  \int_{0}^{t} \nu_{l_{k}}(X_s) \; d \ell_s^{l_{k}}(\|X\|)
$$
corresponds to an outward normal reflection with probability $\alpha_k$ and an inward normal reflection with probability $1-\alpha_{k}$ when $X_t$ hits  $\partial B_{l_k}$ 
(cf. \cite{Tr3}). 
Analogous interpretations hold for the other reflection terms. Thus $\partial B_{l_k}$, $\partial B_{r_k}$ and $\partial B_{m_0}$ can be seen as boundaries where a skew reflection takes place, 
or alternatively as permeable membranes.
\end{remark}

\section{INTEGRATION BY PARTS FORMULA}
\begin{prop}\label{ibp}
Suppose (H1)-(H3) hold. The following integration by parts formula holds for $f,g \in C_0^{\infty}(\R^d)$
\begin{eqnarray*}
-\E(f,g) & = & \int_{\R^d} \left ( \frac{1}{2}\Delta f  + \nabla f \, \cdot \frac{\nabla \rho  }{2\rho}   \right ) g \, \rho \, \phi \, dx 
+ \frac{ \overline{\gamma} -  \gamma}{2}\int_{\partial B_{m_0}}  \nabla  f   \cdot \nu_{m_0} \,   g  \, \rho \, d \sigma_{m_0}  \\
&& \hspace*{-1.5cm}+ \sum_{k \in \mathbb{Z}}\left (  \frac{\gamma_{k+1} - \gamma_{k}}{2} \int_{\partial B_{l_{k}}}   \,\nabla  f   \cdot \nu_{l_{k}} \,   g  \, \rho \, d \sigma_{l_{k}}
+ \frac{\overline{\gamma}_{k+1} - \overline{\gamma}_{k}}{2} \int_{\partial B_{r_{k}}}   \, \nabla  f   \cdot \nu_{r_{k}} \,   g  \, \rho \, d \sigma_{r_{k}}\right ).
\end{eqnarray*}
The last summation is in particular only over finitely many $r_k$, $k\ge 1$, since $f$ has compact support.
\end{prop}
\pf
For $f,g \in C_0^{\infty}(\R^d)$
\begin{eqnarray*}
&&\E(f,g) = \frac{1} {2} \sum_{j=1}^{d}   \int_{\R^d} \partial_j f \partial_j g \;  \left( \rho \, \phi \right) \, dx\\
&=& - \frac{1}{2} \sum_{j=1}^{d} \sum_{k \in \mathbb{Z}} \left(  \, \gamma_{k} \, \int_{A_{k}} \left( \partial_{jj} f  + \partial_j f \, \frac{\partial_j   \rho \,  }{\rho}   \right) g \, \rho \, dx  \, + \, \overline{\gamma}_{k} \, \int_{\hat{A}_{k}} \left( \partial_{jj} f  + \partial_j f \, \frac{\partial_j  \rho  }{\rho}   \right) g \, \rho \, dx \right) \\
&&+ \frac{1} {2} \sum_{j=1}^{d} \sum_{k \in \mathbb{Z}} \left( \int_{A_{k}} \, \gamma_{k} \, \partial_j  \left( \partial_j f \, g  \;  \rho  \right) \, dx
+  \int_{\hat{A}_{k}} \, \overline{\gamma}_{k}\, \partial_j \left( \partial_j f \, g   \;  \rho \right) \, dx \right).\\
\end{eqnarray*}
The first term equals
\begin{equation*}
- \frac{1}{2}   \,  \int_{\R^d} \left( \Delta f  + \nabla f \, \cdot \frac{\nabla \rho  }{\rho}   \right) g \, \rho \, \phi \, dx,
\end{equation*}
and the second term equals
\begin{eqnarray*}
&&\frac{1} {2}  \sum_{j=1}^{d} \sum_{k \in \mathbb{Z}}   \left( \int_{\partial B_{l_{k}}}  \gamma_{k} \, \partial_j  f   \, \nu_{l_{k}}^{j} \,   g  \, \rho \, d \sigma_{l_{k}}
-  \int_{\partial B_{l_{k-1}}}  \gamma_{k} \, \partial_j  f   \, \nu_{l_{k-1}}^{j} \,   g  \, \rho \, d \sigma_{l_{k-1}}
\right)\\
&&+\frac{1} {2}  \sum_{j=1}^{d} \sum_{k \in \mathbb{Z}} \left(  \int_{\partial B_{r_{k}}}  \overline{\gamma}_{k} \, \partial_j  f   \, \nu_{r_{k}}^{j} \,   g  \, \rho \, d \sigma_{r_{k}}
-  \int_{\partial B_{r_{k-1}}}  \overline{\gamma}_{k} \, \partial_j  f   \, \nu_{r_{k-1}}^{j} \,   g  \, \rho \, d \sigma_{r_{k-1}}
\right)
\end{eqnarray*}
\begin{eqnarray*}
&=& -\frac{1} {2} \left( \lim_{k \rightarrow -\infty} \int_{\partial B_{l_{k-1}}}  \gamma_{k} \, \nabla  f   \cdot \nu_{l_{k-1}} \,   g  \, \rho \, d \sigma_{l_{k-1}}
+\sum_{k \in \mathbb{Z}}   \int_{\partial B_{l_{k}}}  ( \gamma_{k+1} - \gamma_{k}) \,\nabla  f   \cdot \nu_{l_{k}} \,   g  \, \rho \, d \sigma_{l_{k}} \right.\\
&&- \left. \lim_{k \rightarrow \infty}  \int_{\partial B_{l_{k+1}}}  \gamma_{k+1} \, \nabla  f   \cdot \nu_{l_{k+1}} \,   g  \, \rho \, d \sigma_{l_{k+1}} \right)\\
&&-\frac{1} {2} \left(  \lim_{k \rightarrow -\infty}  \int_{\partial B_{r_{k-1}}}  \overline{\gamma}_{k} \,\nabla  f   \cdot \nu_{r_{k-1}} \,   g  \, \rho \, d \sigma_{r_{k-1}}
+ \sum_{k \in \mathbb{Z}}   \int_{\partial B_{r_{k}}}  ( \overline{\gamma}_{k+1} - \overline{\gamma}_{k}) \, \nabla  f   \cdot \nu_{r_{k}} \,   g  \, \rho \, d \sigma_{r_{k}} \right)\\
&=&- \frac{1} {2} \left( \sum_{k \in \mathbb{Z}}   \int_{\partial B_{l_{k}}}  ( \gamma_{k+1} - \gamma_{k}) \,\nabla  f   \cdot \nu_{l_{k}} \,   g  \, \rho \, d \sigma_{l_{k}}
- \int_{\partial  B_{m_0}}  \gamma \, \nabla  f   \cdot \nu_{m_0} \,   g  \, \rho \, d \sigma_{m_0}  \right)\\
&&- \frac{1} {2} \left(\int_{\partial B_{m_0}}  \overline{\gamma} \,  \nabla  f   \cdot \nu_{m_0} \,   g  \, \rho \, d \sigma_{m_0}
+ \sum_{k \in \mathbb{Z}}   \int_{\partial B_{r_{k}}}  ( \overline{\gamma}_{k+1} - \overline{\gamma}_{k}) \, \nabla  f   \cdot \nu_{r_{k}} \,   g  \, \rho \, d \sigma_{r_{k}} \right ),
\end{eqnarray*}
because
\begin{equation*}
\lim_{k \rightarrow -\infty} \int_{\partial B_{l_{k-1}}}  \gamma_{k} \, \nabla  f   \cdot \nu_{l_{k-1}} \,   g  \, \rho \, d \sigma_{l_{k-1}}
=\lim_{k \rightarrow -\infty} \sum_{j=1}^{d}   \int_{ B_{l_{k-1}}}  \gamma_{k} \, \partial_j( \partial_j  f  \,   g  \, \rho ) \, dx = 0,
\end{equation*}
by (H1) and Lebesgue, since $\partial_j( \partial_j  f  \,   g  \, \rho ) \in L^1_{loc}(\R^d)$. Similarly
\begin{eqnarray*}
&&\lim_{k \rightarrow \infty}  \int_{\partial B_{l_{k+1}}}  \gamma_{k+1} \, \nabla  f   \cdot \nu_{l_{k+1}} \,   g  \, \rho \, d \sigma_{l_{k+1}}
=\lim_{k \rightarrow \infty} \sum_{j=1}^{d} \int_{ B_{l_{k+1}}}  \gamma_{k+1} \, \partial_j( \partial_j  f  \,   g  \, \rho ) \, dx\\
&&=\sum_{j=1}^{d} \int_{ B_{m_0}}   \gamma \, \partial_j( \partial_j  f  \,   g  \, \rho ) \, dx
= \int_{\partial  B_{m_0}}  \gamma \, \nabla  f   \cdot \nu_{m_0} \,   g  \, \rho \, d \sigma_{m_0},
\end{eqnarray*}
and 
\begin{equation*}
\lim_{k \rightarrow -\infty}  \int_{\partial B_{r_{k-1}}}  \overline{\gamma}_{k} \,  \nabla  f   \cdot \nu_{r_{k-1}} \,   g  \, \rho \, d \sigma_{r_{k-1}}
=  \int_{\partial B_{m_{0}}}  \overline{\gamma} \;  \nabla  f   \cdot \nu_{m_0} \,   g  \, \rho \, d \sigma_{m_0}.
\end{equation*}
\qed
\begin{remark}\label{ibploc}
The integration by parts formula in Proposition \ref{ibp} extends to $f(x)=|\| x \| -a|, a\in \R$, and to the coordinate projections.
\end{remark}
\section{STRICT FUKUSHIMA DECOMPOSITION}\label{4}

A positive Radon measure $\mu$ on $\R^d$ is said to be of {\it finite energy integral} if
\[
\int_{\R^{d}} |f(x)|\, \mu (dx) \le C \sqrt{ \E_1(f,f)} , \; f \in D(\E) \cap C_0(\R^{d}),
\]
where $C$ is some constant independent of $f$, and $C_0(\R^{d})$ is the set of compactly supported continuous functions on $\R^{d}$. A positive Radon measure $\mu$ on $\R^{d}$ is of finite energy integral if and only if there exists a unique function $U_{1} \, \mu\in D(\E)$ such that
\[
\E_{1}(U_{1} \, \mu, f) = \int_{\R^{d}} f(x) \, \mu(dx),
\]
for all $f \in D(\E) \cap C_0(\R^{d})$. $U_{1} \, \mu$ is called $1$-potential of $\mu$. In particular, $R_{1} \mu$ is a $\rho\phi \, dx$-version of $U_{1} \mu$. The measures of finite energy integral are denoted by $S_0$.\\
Let further
\[
S_{00} : = \{\mu\in S_0 \,|\;  \mu(\R^{d})<\infty, \|U_{1} \mu\|_{\infty}<\infty \},
\]
where $\| f \|_{\infty} : = \inf\{c>0\,|\; \int_{\R^d}  1_{\{ \,|f|>c \, \} } \, \rho \phi \, dx = 0   \}$ and define for $l \in (0, \infty)$
\[
C_{0}^{\infty}( \overline{B}_{l}) : = \{\, f : \overline{B}_{l} \rightarrow \R \,|\;  \exists g \in C_{0}^{\infty}(\R^d) \; with \;  f=g  \; on \; \overline{B}_{l}   \}.
\]

\begin{lemma}\label{L;LEM2}
Suppose (H1)-(H3) hold. Then for  $l \in (0, \infty)$ and $f \in C_{0}^{\infty}( \overline{B}_{l})$,
\[
\int_{\partial B_l} |\,  f \,| \, \rho \, d\sigma_{l} \le C(l) \, \left( \int_{B_{l}} \|\nabla \,f\,\|^2 \, \rho \, \phi  \, dx + \int_{B_{l}} |\,f\,|^2 \, \rho \, \phi \, dx \right)^{1/2},
\]
where $C : (0, \infty) \rightarrow \R$ is an increasing function.
In particular, for any $f \in D(\E)$
\[
\int_{\partial B_l} |\,  f \,| \, \rho \, d\sigma_{l} \le C(l) \, \| \,f\, \|_{D(\E)}.
\]
\end{lemma}
\pf
Since $\partial B_{l}$ has Lipschitz boundary (actually $C^{\infty}$-boundary), we can see from \cite[Theorem 1 in section 4.3]{EvGa}  (and  \cite[Theorem 5.1 (i)]{Tr1}) that there exists a constant $\sqrt{2}$ independent of $l$, such that for $f \in C_{0}^{\infty}( \overline{B}_{l})$
\[
\int_{\partial B_l} |\,  f \,| \, \rho \, d\sigma_{l}  \le \sqrt{2} \int_{B_{l}} \big( \|\, \nabla f\,\| \, \rho \, + \, 2 |\, \xi \, f  \,| \,  \| \nabla \xi \| + |\,f\,| \, \rho  \big)\, dx
\]
\[
\le \sqrt{2}  \left[ \left( \int_{B_{l}} \|\nabla \,f\,\|^2 \, \rho \, dx  \right)^{1/2} \  \big(\rho \, dx(B_{l}) \big)^{1/2} + 2 \left( \int_{B_{l}} |\,f\,|^2 \, \rho \, dx  \right)^{1/2}  \,
 \left( \int_{B_{l}} \|\, \nabla \xi\,\|^2  \, dx  \right)^{1/2} \right.
\]
\[
\left. + \left( \int_{B_{l}} |\,f\,|^2 \, \rho \, dx  \right)^{1/2} \  \big(\rho \, dx(B_{l}) \big)^{1/2} \right]
\]
\[
\le \sqrt{\frac{8}{\delta_{l}}} \left(  \  \big(\rho \, dx(B_{l}) \big)^{1/2} + || \nabla \xi  ||_{L^{2}(B_{l}, dx)}   \right) \,  \left( \int_{B_{l}} \|\nabla \,f\,\|^2 \, \rho \, \phi  \, dx + \int_{B_{l}} |\,f\,|^2 \, \rho \, \phi \, dx \right)^{1/2}.
\]
This is the first statement. Let $f \in C_{0}^{\infty}(\R^d)$. Then, by the above, since $f $ restricted to $\overline{B}_{l}$ is in $C_{0}^{\infty}(\overline{B}_{l})$,
\[
\int_{\partial B_l} |\,  f \,| \, \rho \, d\sigma_{l} \le  C(l) \, \left( \int_{B_{l}} \|\nabla \,f\,\|^2 \, \rho \, \phi  \, dx + \int_{B_{l}} |\,f\,|^2 \, \rho \, \phi \, dx \right)^{1/2}
\le C(l) \, ||f ||_{D(\E)}.
\]
Since $C_{0}^{\infty}(\R^{d})$ is dense in $D(\E)$, the second statement follows. 
\qed

A positive Borel measure $\mu$ on $\R^d$ is said to be smooth in the strict sense if there exists a sequence $(E_k)_{k\ge 1}$ of Borel sets increasing to $\R^d$ such that $1_{E_{k}} \cdot \mu \in S_{00}$ for each $k$ and
\[
\P_{x} ( \lim_{k \rightarrow \infty} \sigma_{ \R^d \setminus E_{k} }  \ge \zeta ) =1 \;, \;\; \forall x \in \R^d.
\]
Here  $\sigma_{B}:=\inf\{t>0\,|\, X_t\in B\}$ for $B \in \mathcal{B}(\R^{d})$. The totality of the smooth measures in the strict sense is denoted by $S_{1}$ (see \cite{FOT}).
\begin{lemma}\label{L;LEM3}
Suppose (H1)-(H4) hold. Let $l  \in (0, \infty)$. Then, for any relatively compact open set $G$,  $1_{G} \cdot \rho d \sigma_{l} \in S_{00}$. In particular, $\rho d \sigma_{l} \in S_{1}.$

\end{lemma}

\pf
Let $l  \in (0, \infty)$. By Lemma $\ref{L;LEM2}$, $\rho \, d\sigma_{l} \in S_{0}$. Let us first show that $\rho \,d \sigma_{l} \in S_{1}$ with respect to an increasing sequence of open sets $(E_{k})_{k \ge 1}$. Choose $\varphi, \, \overline{\varphi} \in L_{b}^{1}(\R^d, \rho \phi \, dx)$, $0 < \varphi, \, \overline{\varphi} \le 1$ $\rho \phi \, dx$-a.e. By assumption (H4) $R_1(\rho d\sigma_{l})$ is continuous. Since furthermore $R_1 \varphi$ is continuous and $R_1 \varphi$ is strictly positive, it follows that
\[
E_k : = \{\, R_1(\rho d\sigma_{l}) < k^2 R_1 \varphi  \,  \} ,\;\; k \ge 1,
\]
are open sets that increase to $\R^d$. Choosing the constant function $1 \in C_{0}^{\infty}(\overline{B}_{l})$ in Lemma $\ref{L;LEM2}$ we see that $\rho d\sigma_{l}$ is finite. Then, clearly $1_{E_k} \cdot \rho d\sigma_{l} $ is also finite for all $k \ge 1$. So, it remains to show that the corresponding 1-potentials $U_1(1_{E_k} \cdot \rho d\sigma_{l})$ are $\rho \phi \, dx$-essentially bounded. Let $\left( G_{1} \overline{\varphi}  \right)_{E_k}$ be the 1-reduced function of $G_1 \overline{\varphi}$ on $E_k$ as defined in \cite{Tr2}. Then $\EE_{\cdot} \left[ \int_{\sigma_{E_k}}^{\infty} e^{-t} \overline{\varphi}(X_t) \, dt  \right]$ is a $\rho\phi\, dx$-version of $\left( G_{1} \overline{\varphi}  \right)_{E_k}$. We have (for intermediate steps see \cite[p.416]{Tr2})
\begin{eqnarray*}
&&\int_{\R^d}\overline{\varphi} U_1 \left( 1_{E_k} \cdot \rho d\sigma_{l}    \right)   \rho\phi\, dx
= \int_{\partial B_{l}} R_{1} \overline{\varphi} \;\; 1_{E_k} \cdot \rho d\sigma_{l}\\
&&= \int_{\partial B_{l}} \EE_{\cdot} \left[ \int_{\sigma_{E_k}}^{\infty} e^{-t} \overline{\varphi}(X_t) \, dt  \right]    1_{E_k} \cdot \rho d\sigma_{l}
= \E_{1} \left( \left( G_{1} \overline{\varphi}  \right)_{E_k} ,  U_1 ( 1_{E_k} \cdot \rho d\sigma_{l})  \right)\\
&&= \E_{1} \left( \left( G_{1} \overline{\varphi}  \right)_{E_k} ,  U_1 ( 1_{E_k} \cdot \rho d\sigma_{l}) \wedge k^2 G_{1} \varphi  \right)\\
&&\le  \E_{1} \left( G_1 \overline{\varphi} ,  U_1 ( 1_{E_k} \cdot \rho d\sigma_{l}) \wedge k^2 G_{1} \varphi  \right)
= \int_{\R^d}\overline{\varphi}  \left ( U_1 ( 1_{E_k} \cdot \rho d\sigma_{l}) \wedge k^2 G_{1} \varphi \right ) \,\rho\phi\, dx.
\end{eqnarray*}
This implies that $U_1 ( 1_{E_k} \cdot \rho d\sigma_{l}) \le k^2 \;\; \rho \phi \, dx$-a.e. Hence, $ 1_{E_k} \cdot \rho d\sigma_{l} \in S_{00}$ for all $k \ge 1$. Since moreover $\P_{x} ( \lim_{k \rightarrow \infty} \sigma_{ \R^d \setminus E_k }  \ge \zeta ) =1 , \; \forall x \in \R^d$ is easily deduced from the absolute continuity condition,
we finally obtain $\rho d\sigma_{l} \in S_{1}$ with respect to $(E_{k})_{k \ge 1}$.\\
For a relatively compact open set $G$, we know that there exists $k_{0} \in \N$ with $G \subset \overline{G} \subset E_{k_{0}}$. Hence, 
 $U_1 \left(1_{G} \cdot \rho \, d\sigma_{l} \right) \le U_1 (1_{E_{k_{0}}} \cdot \rho \, d\sigma_{l} ) \le k_{0}^{2} \, G_{1} \varphi \le k_{0}^{2}$.
Therefore, $1_{G} \cdot \rho d \sigma_{l} \in S_{00}$.
\qed

By Lemma $\ref{L;LEM3}$, we know that $\rho d\sigma_{r} \in S_{1}$ for any $r \in (0, \infty)$. Hence, by \cite[Theorem 5.1.7]{FOT} there exists a unique $(\bar{\ell}^{r}_t)_{t \ge 0} \in A_{c,1}^{+}$ with Revuz measure $\rho d \sigma_{r}$, such that
\[
\EE_x \left[ \int_{0}^{\infty} e^{-t} \, d \bar{\ell}_{t}^{r}   \right] = R_1(\rho d\sigma_{r})(x) \, , \;\; \forall x \in \R^d.
\]
Here, $A_{c,1}^{+}$ denotes the positive continuous additive functionals in the strict sense. 
\begin{thm}\label{T;SUR}
Suppose (H1)-(H4) hold. Then, for any relatively compact open set $G$, $1_{G} \cdot \mu \in S_{00}-S_{00}$, where
\[
\mu=\sum_{k \in \mathbb{Z}} \left (  \frac{ \gamma_{k+1} - \gamma_{k}}{2} \,   \rho \,d \sigma_{l_{k}} + \frac{ \overline{\gamma}_{k+1} - \overline{\gamma}_{k} }{2} \, \rho \,d 
 \sigma_{r_{k}}\right ) +  \frac{ \overline{\gamma} -  \gamma }{2} \, \rho \, d \sigma_{m_0}.
\]
In particular $\mu \in S_{1}-S_{1}$.
\end{thm}
\pf
It suffices to show that $\mu_{n} \in S_{00}$ for any $n \in \N, \; n > m_0$, where
\[
\mu_{n} : =1_{G} \cdot \left( \sum_{k \in \mathbb{Z}} |\,  \gamma_{k+1} - \gamma_{k}  \, | \, \rho \,d \sigma_{l_{k}} + \sum_{ \{ k \in \mathbb{Z} :  r_{k}< n  \}}  |\,  \overline{\gamma}_{k+1} - \overline{\gamma}_{k} \,| \, \rho \,d \sigma_{r_{k}} + | \,  \overline{\gamma} -  \gamma \, | \, \rho \, d \sigma_{m_0} \right). 
\]
First we show that $\mu_{n} \in S_{0}$. Let $f \in C_{0}^{\infty}(\R^d)$.
Then, by Lemma $\ref{L;LEM2}$
\[
\int_{\R^{d}}  |\,f\,|    \,d \mu_{n}
\le C(n) \,  \left(  \sum_{k \in \mathbb{Z}} |\,  \gamma_{k+1} - \gamma_{k}  \, | +  \sum_{ \{ k \in \mathbb{Z} :  r_{k}< n  \}}  |\,  \overline{\gamma}_{k+1} - \overline{\gamma}_{k} \,| + | \, \overline{\gamma} -  \gamma \, |   \right)\; ||f||_{D(\E)},
\]
where $C(n)$ is as in Lemma $\ref{L;LEM2}$.
Now, we show $\mu_{n} \in S_{00}$.
Let $f \in C_{0}^{\infty}(B_{n})$ such that $f=1$ on $\overline{B}_{{n-\epsilon}}$ where $\epsilon>0$ is small enough to satisfy $( n-\epsilon) > max_{ \{ k \in \mathbb{Z} :  r_{k}< n  \}} r_{k}$. Then, by Lemma $\ref{L;LEM2}$ again
\[
\mu_{n}( \R^{d} )
\le C(n) \left(  \sum_{k \in \mathbb{Z}} |\,  \gamma_{k+1} - \gamma_{k}  \, | +  \sum_{ \{ k \in \mathbb{Z} :  r_{k}< n  \}}  |\,  \overline{\gamma}_{k+1} - \overline{\gamma}_{k} \,| + | \, \overline{\gamma} -  \gamma \, |   \right)  \left( \int_{B_{n}} \, \rho \phi \,dx \right)^{1/2}< \infty.
\]
By the proof of Lemma $\ref{L;LEM3}$, we can see that 
\[
U_{1} \left(1_{E_{k}} \cdot \rho \, d\sigma_{l} \right) \le k^{2}, \;\; k \ge 1,
\]
independently of $l \in (0, \infty)$. For any relatively compact open set $G$, there exists $k_{0} \in \N$ such that $G \subset \overline{G} \subset E_{k_{0}}$.
Since $U_1 \left(1_{G} \cdot \rho \, d\sigma_{l} \right) \le U_1 (1_{E_{k_{0}}} \cdot \rho \, d\sigma_{l} )$ for any $l$,
we obtain  for $\rho\phi \,dx$-a.e. $x \in \R^d$
\begin{eqnarray*}
U_1\mu_{n}(x) &\le& \sum_{k \in \mathbb{Z}} |\,  \gamma_{k+1} - \gamma_{k}  \, |  \; U_1 \left(1_{G} \cdot  \rho \,  d\sigma_{l_{k}} \right) (x) \\
&+& \sum_{ \{ k \in \mathbb{Z} :  r_{k}< n  \} } |\,  \overline{\gamma}_{k+1} - \overline{\gamma}_{k}  \, |  \; U_1 \left(1_{G} \cdot \rho \,  d\sigma_{r_{k}} \right) (x) 
+ |\,  \overline{\gamma} - \gamma  \, |   \; U_1 \left(1_{G} \cdot \rho \,  d\sigma_{m_{0}} \right) (x)  \\
&\le& k^2_{0} \;  \left(  \sum_{k \in \mathbb{Z}} |\,  \gamma_{k+1} - \gamma_{k}  \, | +  \sum_{ \{ k \in \mathbb{Z} :  r_{k}< n  \} } |\,  \overline{\gamma}_{k+1} - \overline{\gamma}_{k}  \, |  +    |\,  \overline{\gamma} - \gamma  \, |    \right) < \infty.
\end{eqnarray*}
Therefore, $\mu \in S_{00}$.
\qed

\begin{remark}\label{uniquenessofsum}
Let $E_k, k\ge 1$, be open sets such that $E_k\nearrow \mathbb{R}^d$ and let $\mu=\mu_A, \mu_n=\mu_{A^n}\in S_1$ w.r.t. $(E_k)_{k\ge 1}$, $A, A^n \in A_{c,1}^{+}$, $n\ge 1$. 
If $\mu_A=\sum_{n\ge 1}\mu_{A^n}$, then 
$A=\sum_{n\ge 1} A^n$, since $R_1 (f d\mu_A)(x)=\sum_{n\ge 1}R_1(f d\mu_{A^n})(x)$ for any $x\in \mathbb{R}^d$, $f\in C_0(\mathbb{R}^d)$. \\
\end{remark}

\begin{thm}\label{C;SM1}
Suppose (H1)-(H4) hold. For any relatively compact open set $G$ and $j=1, \dots , d$, 
\[
1_{G} \cdot \partial_j \rho \, \phi \, dx \in S_{00} - S_{00}.
\]
In particular $\partial_j \rho \, \phi \, dx \in S_{1} - S_{1}$, $j=1, \dots , d$.
\end{thm}
\proof
It suffices to show that $1_{G} \cdot  | \partial_j \rho | \, \phi \, dx \in S_{00}$.
By Remark $\ref{rem2}$, it is easy to see that  $1_{G} \cdot  | \partial_j \rho | \, \phi \, dx \in S_{0}$ and that $1_{G} \cdot  | \partial_j \rho | \, \phi \, dx (\R^d) < \infty$.
We can show that 
\[
\left\| U_1(1_{G} \cdot  | \partial_j \rho | \, \phi \, dx ) \right\|_{\infty} < \infty
\]
by replacing the $E_{k}$ in the proof of Lemma $\ref{L;LEM3}$ with $E_{k}^{\prime} := \{\, R_1(1_{G} \cdot  | \partial_j \rho | \, \phi \, dx) < k^2 R_1 \varphi  \,  \}$.
\qed

\section{PROOF OF THE MAIN THEOREM}\label{5}

\text{\it Proof of Theorem \ref{T;MAIN}}.
(i) Applying \cite[Theorem 5.5.5]{FOT} to $(\E,D(\E))$ and to the coordinate projections which are in $D(\E)_{b,loc}$, the identification of the martingale part as Brownian motion is 
easy. Concerning the drift part the strict decomposition holds true by Lemma \ref{L;LEM3}, Remark \ref{ibploc},  Theorem $\ref{T;SUR}$, Remark \ref{uniquenessofsum} and Theorem $\ref{C;SM1}$. Note that equation (\ref{SDE}) holds for all $t\ge 0$ by (\ref{HPcons}). \\
(ii) Let $f(x):=\|x\|$, $x\in\mathbb{R}^d$. Then $\partial_j f$ is everywhere bounded by one (except in zero). Thus applying \cite[Theorem 5.5.5]{FOT} to $f$, which is in $D(\E)_{b,loc}$, we obtain similarly to (i)
\begin{equation}\label{SDEf}
\|X_t\| = \|x\| + B_t + \int^{t}_{0} \frac{X_s}{\|X_s\|}\cdot\frac{\grad{\rho}}{2\rho} (X_s) \, ds +  \overline{N}_{t}, 
\end{equation}
$ \P_x $ -a.s. for any $x \in \R^d$,  $t\ge 0$,  where $B$ is a standard one dimensional Brownian motion starting from zero and
\begin{equation*}
\overline{N}_{t} : =  \sum_{k \in \mathbb{Z}} \left (\frac{ \gamma_{k+1} - \gamma_{k} } { \gamma_{k+1} + \gamma_{k} }  \ \ell_t^{l_{k}}
+ \frac{ \overline{\gamma}_{k+1} - \overline{\gamma}_{k} } { \overline{\gamma}_{k+1} + \overline{\gamma}_{k} } \ \ell_t^{r_{k}}\right )
+ \frac{ \overline{\gamma}  - \gamma }{ \overline{\gamma}  + \gamma } \ \ell_t^{m_0}.
\end{equation*}
Therefore, the first statement follows.
In particular, we may apply the symmetric It\^o-Tanaka formula (see \cite[VI. (1.25)]{RYor}) and obtain
\begin{equation}\label{Tanaka}
\big | \|X_t\|-a\big | = \big |\| x \|-a\big | + \int_{0}^t sign(\|X_s\|-a)d\|X_s\| + \ell_t^a(\|X\|),
\end{equation}
$ \P_x $ -a.s. for any $x \in \R^d$,  $t\ge 0$, where $sign$ is the point symmetric sign function. 
Let $h(x):=\big | \|x\|-a \big |$, $a\in\{m_0, l_k,r_k : k\in \mathbb{Z}\}, x\in\mathbb{R}^d$. Then $\partial_j h$ is everywhere bounded by one (except in $a$).
Thus, applying \cite[Theorem 5.5.5]{FOT} to $h$, which is in $D(\E)_{b,loc}$, we obtain again similarly to (i)
\begin{equation}\label{TanakaDF}
\big | \|X_t\|-a\big | = \big |\| x \|-a\big | + \int_{0}^t sign(\|X_s\|-a)d\|X_s\| + \ell_t^a,
\end{equation}
$ \P_x $ -a.s. for any $x \in \R^d$, $t\ge 0$. Comparing (\ref{Tanaka}) and (\ref{TanakaDF}), we get the result.

\qed

\begin{remark}\label{R;CONS1}
(see \cite[Theorem 4.7.1 (i), (iii), and Exercise 4.7.1]{FOT}) If $(\E,D(\E))$ is irreducible, then for any nearly Borel non-exceptional set $B$,
\[
\P_{x} ( \sigma_{B} < \infty) >0, \;\;  \forall x \in \R^d.
\]
If $(\E,D(\E))$ is irreducible and recurrent, then for any nearly Borel non-exceptional set $B$,
\[
\P_{x} ( \sigma_{B}  \circ \theta_{n} < \infty, \forall n \ge 0 ) =1, \;\;  \forall x \in \R^d.
\]
Here $(\theta_{t})_{t\ge 0}$ is the shift operator. Moreover, in this case any excessive function is constant. In particular, the ergodic Theorem \cite[Theorem 4.7.3 (iv)]{FOT} holds. A sufficient condition for recurrence is given by
\[
\int_1^{\infty} \frac{r}{ \rho \phi \,dx (B_{r}) }dr= \infty,
\]
see \cite[Theorem 3]{St1}.

\end{remark}
{\bf Acknowledgments.}  
The research of Jiyong Shin and of the corresponding author Gerald Trutnau was supported by Basic Science Research Program through the National Research Foundation of Korea(NRF) funded by the Ministry of Education, Science and Technology (MEST) (2013028029) and Seoul National University Research Grant 0450-20110022.

\addcontentsline{toc}{chapter}{Index}
\bigskip
\noindent
{\it Department of Mathematical Sciences and Research Institute of Mathematics of Seoul National University,
599 Gwanak-Ro, Gwanak-Gu, Seoul 151-747, South Korea}

\end{document}